\documentclass[12pt]{amsart}
\usepackage{amsthm,amssymb,amscd,amsmath,latexsym,indentfirst,color,amsfonts}
\usepackage[mathcal]{eucal}
\usepackage{geometry,amsthm,graphics,tabularx,shapepar,float}
\usepackage[all]{xypic}
\usepackage{tikz}
\usepackage{caption}
\usetikzlibrary{calc}
\theoremstyle{plain}
\newtheorem{theorem}{Theorem}[section]
\newtheorem{prop}[theorem]{Proposition}
\newtheorem{coro}[theorem]{Corollary}
\newtheorem{lemma}[theorem]{Lemma}
\newtheorem{ex}[theorem]{Example}

\newtheorem{rem}[theorem]{Remark}

\newcommand{\ble}{\begin {lemma}}
\newcommand{\ele}{\end {lemma}}
\newcommand{\bprp}{\begin {Proposition}}
\newcommand{\eprp}{\end {Proposition}}
\newcommand{\bthm}{\begin {theorem}}
\newcommand{\ethm}{\end {theorem}}
\newcommand{\bco}{\begin {coro}}
\newcommand{\eco}{\end {coro}}
\newcommand{\bex}{\begin {ex}}
\newcommand{\eex}{\end {ex}}

\newcommand{\be}{\begin {equation}}
\newcommand{\ee}{\end {equation}}
\newcommand{\bp}{\begin {proof}}
\newcommand{\ep}{\end {proof}}
\newcommand{\bee}{\begin {equation*}}
\newcommand{\eee}{\end {equation*}}
\newcommand{\rt}{\rightarrow}
\newcommand{\lb}{\label}

\newcommand{\al}{\alpha}

\newcommand{\la}{\lambda}

\newcommand{\pt}{\partial}

\begin{document}

\title{Some progress in the Dixmier Conjecture}
\author [G. Han] {Gang Han}
\address [G. Han]{School of Mathematics   \\ Zhejiang  University \\ Hangzhou, 310027, China} \email{mathhgg@zju.edu.cn}
 \author[B. Tan]{Bowen Tan }
  \address [B. Tan]{School of Mathematics \\ Zhejiang  University\\Hangzhou, 310027, China }\email{supertomato696@gmail.com }
\date{Sept. 28, 2022 }
\subjclass[2010]{16S32, 16W20}
\keywords{ Weyl algebra, Dixmier Conjecture, Newton polygon}
\begin{abstract} Let $p$ and $q$, where $pq-qp=1$, be the standard generators of the first Weyl
algebra  $A_1$ over a field of characteristic zero. Then the spectrum of the inner derivation $ad(pq)$ on $A_1$ are exactly the set of integers. The algebra $A_1$ is a $\mathbb{Z}$-graded algebra with each $i$-component being the $i$-eigenspace of $ad(pq)$, where $i\in \mathbb{Z}$. The Dixmier Conjecture says that  if some elements $z$ and $w$ of $A_1$ satisfy $zw-wz=1$,  then they generate $A_1$.  We show that if either $z$ or $w$ possesses no component belonging to the negative spectrum of $ad(pq)$, then the Dixmier Conjecture holds. We give some generalization of this result, and some other useful criterions for $z$ and $w$ to generate $A_1$. An important tool in our proof is the Newton polygon for elements in $A_1$.
\end{abstract}
\thanks{Bowen Tan is the corresponding author.}
 \maketitle

\section{Introduction}
 \setcounter{equation}{0}\setcounter{theorem}{0}

Let $K$ be a field of Characteristic 0. Let $A_1$ be the first Weyl algebra over $K$ generated by $p,q$ with \[[p,q]=pq-qp=1.\]

In the influential paper \cite{d}, Dixmier studied the structure of $A_1$ in great detail, and in the end he posed six problems. Joseph gave some characterization  for semisimple and nilpotent elements of $A_1$ and solved two problems in \cite{jo}.

 The first problem of \cite{d} asks if
every endomorphism of the Weyl algebra $A_1$ is an automorphism, i.e., given elements
$z,w$ of $A_1$ such that $[z,w] = 1$, do they generate the algebra $A_1$?  A similar problem
for the $n$-th Weyl algebra is called the Dixmier Conjecture. It is open for all $n\ge 1$. In this paper, the Dixmier Conjecture for $A_1$ will be abbreviated by DC.

People have made many explorations in DC.

In \cite{mv}, the authors showed that DC holds for a particular type of endomorphism of $A_1$. By taking the same idea further, Moskowicz gave some equivalent formulation of DC \cite{mo}.

The spectrum of the inner derivation $ad(pq)$ on $A_1$ are exactly the set of integers. $A_1$ is a $\mathbb{Z}$-graded algebra with each $i$-component $D_i$ being the $i$-eigenspace of $ad(pq)$, $i\in \mathbb{Z}$. An element $z\in D_i$ is said to be homogeneous of degree $i$. One has \[A_1=\oplus_{i\in \mathbb{Z}}D_i,\ and\ D_i\cdot D_j\subseteq D_{i+j}, \forall i,j\in \mathbb{Z}.\] It follows that $D_{\ge 0}=\oplus_{i\ge 0}D_i$ is a subalgebra of $A_1$.

In \cite{b}, Bavula showed that the eigenvector algebra of $zw$, i.e., the subalgebra of $A_1$ generated by all the eigenvectors of $ad(zw)$ in $A_1$, is equal to the subalgebra $K\langle z,w\rangle$ of $A_1$ generated by $z$ and $w$.

In \cite{ggv}, the authors showed that if $z,w$ of $A_1$ satisfying $[z,w] = 1$, then the centralizer $C(z)$ of $z$ in $A_1$ equals $K[z]$.

If $z$ and $w$ do not generate $A_1$, the authors in \cite{vgg} showed that the
greatest common divisor of the total degrees of $z$ and $w$ are $>15$, where $z$ and $w$ are written as polynomials in $p$ and $q$.

In \cite{bl} the authors prove that DC holds if  $z$ and $w$  are sums of no more than two homogeneous elements of $A_1$.

In this paper we show that if $z$ and $w$ satisfy certain conditions then they generate $A_1$. In particular we prove the following  main result in Theorem \ref{33}.
\bthm
Assume that $z,w\in A_1,$ $[z,w]=1$, and either $z$ or $w$ is in $D_{\ge 0}$. Then  $z$ and $w$ generate $A_1$.
\ethm

This result is generalized in Corollary \ref{98}. Proposition \ref{93} and Theorem \ref{56} give some other criterions for  $z$ and $w$ to generate $A_1$. As a corollary,  the main result in \cite{bl} is generalized in Theorem \ref{97}.

Our proof are partially motivated by the idea in \cite{jo} of Joseph. The Newton polygon for any element in $A_1$ is defined and plays important role in the proof. Our definition of Newton polygon is slightly different to that defined in \cite{ml2}. Some important observations of Newton polygons of elements in $A_1$ can be found in \cite{gz}. \\[3mm]

Here are some of the conventions of the paper.

$X,Y$ denote indeterminants;

$i,j,k,l,m,n,r,s$ denote integers; %$\rho,\sigma,\tau$ denote real numbers;

For positive integers $m$ and $n$, $gcd(m,n)$ is the positive greatest common devisor of $m$ and $n$;

%$\mathbb{R}_{\ge 0}=\{a\in \mathbb{R}|a\ge 0\}$;

$K$ is a field of characteristic 0; $\alpha,\beta,\gamma,\lambda,\mu$ denotes elements in $K$; $K^*=K\setminus\{0\}$;

For $z\in A_1$, $C(z)=\{w\in A_1|zw=wz\}$ is the centralizer of $z$ in $A_1$; $K[z]=\{f(z)|f\in K[X]\}$;

$\mathbb{R}^2$ denotes the Euclidean plane with the Cartesian coordinate system $Oxy$;

$(\mathbb{R}_{\ge 0})^2=\{(x,y)\in \mathbb{R}^2|x\ge 0,y\ge 0\}$;

$V_+=\{(x,y)\in (\mathbb{R}_{\ge 0})^2|y\ge x\}, V_-=\{(x,y)\in (\mathbb{R}_{\ge 0})^2|y\le x\}$;

For a nonempty set $S\subseteq \mathbb{R}^2$,\\[1mm]

\centerline{$Convex(S)$ is the convex hull of $S$ in $\mathbb{R}^2$;}
  \[Cone(S)=\{t(x,y)| t\in \mathbb{R}_{\ge 0}, (x,y)\in S\}.\]

\section{Some preliminaries}
 \setcounter{equation}{0}\setcounter{theorem}{0}
\noindent\textbf{2.1.} Recall that $K$ is a field of Characteristic 0. For $f,g\in K[X,Y]$, define their Poisson product
\be \lb{7}
\{f,g\}=\left |
\begin{array}{ll}
{\pt f}/{\pt X} & {\pt f}/{\pt Y} \\
{\pt g}/{\pt X} & {\pt g}/{\pt Y}
\end{array}
\right |.
\ee
In particular, one has
\be \lb{54}
\{X^iY^j,X^k Y^l\}=\left |
\begin{array}{ll}
i & j \\
k & l
\end{array}
\right | X^{i+k-1}Y^{j+l-1}.
\ee

For $(\rho,\sigma)\in \mathbb{R}^2\setminus \{(0,0)\}$ and $\tau\in \mathbb{R}$, if $f\in K[X,Y]$ satisfies \[\rho X\frac{\pt f}{\pt X}+\sigma Y\frac{\pt f}{\pt Y}=\tau f,\] then $f$ is said to be $(\rho,\sigma)$-\textbf{homogeneous of degree} $\tau$, and is denoted by  $\mathbf{v}_{\rho,\sigma}(f)=\tau$. Note that if $f$ is $(\rho,\sigma)$-homogeneous of degree $\tau$, then $f$ is $(\iota\rho,\iota\sigma)$-homogeneous of degree $\iota\tau$, for any $\iota\in \mathbb{R}\setminus \{0\}$. If $(r,s)\in \mathbb{Z}^2$, then for $X^iY^j\in K[X,Y]$, \[\mathbf{v}_{r,s}(X^iY^j)=ri+sj\in \mathbb{Z}.\]

By this definition, 0 is $(\rho,\sigma)$-homogeneous of degree $\tau$, for any $\tau$ in $\mathbb{R}$; a nonzero constant polynomial is always $(\rho,\sigma)$-homogeneous of degree 0.

For $(\rho,\sigma)\in \mathbb{R}^2\setminus \{(0,0)\}\ and\ \tau\in \mathbb{R}$, let
\[ V_{\rho,\sigma}(\tau)=\{f\in K[X,Y]|\rho X\frac{\pt f}{\pt X}+\sigma Y\frac{\pt f}{\pt Y}=\tau f\}. \] Then $V_{\rho,\sigma}(\tau)$ is a linear subspace of $K[X,Y]$.
\begin{prop}\lb{78}(Proposition 2.1 of \cite{jo})
If $f\in V_{\rho,\sigma}(u), g\in V_{\rho,\sigma}(v)$ with $u,v\in \mathbb{R}$, then \[fg\in V_{\rho,\sigma}(u+v),\ and\ \{f,g\}\in V_{\rho,\sigma}(u+v-\rho-\sigma).\]\\[1mm]
\end{prop}

For $f\in K[X,Y]\setminus\{0\}$, assume $f=\sum \alpha_{ij} X^iY^j$. Define \[E(f)=\{(i,j)\in \mathbb{Z}^2| \alpha_{ij}\ne 0\}.\]
Let $Convex(E(f))$ be the convex hull of $E(f)$ in $\mathbb{R}^2$. It is usually a solid polygon.
For $(\rho,\sigma)\in \mathbb{R}^2\setminus \{(0,0)\}$, we also set \[\mathbf{v}_{\rho,\sigma}(f)=sup\{\rho i+\sigma j|(i,j)\in E(f)\},\] which is usually called the $(\rho,\sigma)$-degree of $f$. Note that if $f$ is a nonzero $(\rho,\sigma)$-homogeneous polynomial then $\mathbf{v}_{\rho,\sigma}(f)$ agrees with the previous definition.  Let
\[E(f;\rho,\sigma)=\{(i,j)\in E(f)|\rho i+\sigma j=\mathbf{v}_{\rho,\sigma}(f).  \}\]

Let
\[\mathbf{f}_{\rho,\sigma}(f)=\sum_{(i,j)\in E(f;\rho,\sigma)} \alpha_{ij} X^iY^j. \]

Assume that $f\in K[X,Y]\setminus \{0\}$, and $(\rho,\sigma)\in \mathbb{R}^2\setminus \{(0,0)\}$, then there is a unique decomposition \be\lb{77}f=\sum_{i=0}^k f_i\ee with $f_i\ne 0$ being $(\rho,\sigma)$-homogeneous for all $i$ and $\mathbf{v}_{\rho,\sigma}(f_i)>\mathbf{v}_{\rho,\sigma}(f_{i+1})$ for $i=0,\cdots,k-1$. We will refer to \eqref{77} as the $(\rho,\sigma)$-\textbf{homogeneous decomposition} of $f$, and $\mathbf{f}_{\rho,\sigma}(f)=f_0$ is the $(\rho,\sigma)$-\textbf{leading component}. Note that for $\tau>0$, one has $\mathbf{f}_{\tau\rho,\tau\sigma}(f)=\mathbf{f}_{\rho,\sigma}(f)$.

The set $Convex(E(\mathbf{f}_{\rho,\sigma}(f)))$ is either a vertex or an edge of $Convex(E(f))$.

\begin{prop}\lb{24}
Let $f,g\in K[X,Y]\setminus \{0\}$, $(\rho,\sigma)\in \mathbb{R}^2\setminus \{(0,0)\}$, $\mathbf{f}_{\rho,\sigma}(f)=f_0$, and $\mathbf{f}_{\rho,\sigma}(g)=g_0$. Then

(1) $\mathbf{f}_{\rho,\sigma}(fg)=f_0g_0$;

(2) If $\{f_0,g_0\}\ne 0$ then $\mathbf{f}_{\rho,\sigma}(\{f,g\})=\{f_0,g_0\}$.
\end{prop}
\bp
Assume that  $f=\sum_{i=0}^k f_i$ and  $g=\sum_{j=0}^l g_j$ are the respective $(\rho,\sigma)$-homogeneous decomposition of $f$ and $g$, where
 $\mathbf{v}_{\rho,\sigma}(f_i)=\tau_i$, and $\tau_i>\tau_{i+1}$ for $i=0,\cdots,k-1$; $\mathbf{v}_{\rho,\sigma}(g_j)=\nu_j$, and $\nu_j>\nu_{j+1}$ for all $j=0,\cdots,l-1$.

(1) One has \[fg=f_0 g_0+\sum_{i+j\ge 1}f_i g_j,\]

where $f_ig_j$ are $(\rho,\sigma)$-homogeneous of degree $\tau_i+\nu_j$, for all $i,j$.

As $f_0 g_0\ne 0$, for any $i,j$ with $i+j\ge 1$, \[\tau_0+\nu_0=\mathbf{v}_{\rho,\sigma}(f_0g_0)>\mathbf{v}_{\rho,\sigma}(f_ig_j)=\tau_i+\nu_j,\] one has $\mathbf{f}_{\rho,\sigma}(fg)=f_0g_0$.

(2) One has \[\{f,g\}=\{f_0,g_0\}+\sum_{i+j\ge 1}\{f_i,g_j\}, \] where $\{f_i,g_j\}$ are $(\rho,\sigma)$-homogeneous of degree $\tau_i+\nu_j-(\rho+\sigma)$, for all $i,j$.

By assumption, $\{f_0,g_0\}\ne 0$. Since for any $i,j$ with $i+j\ge 1$, \[\tau_0+\nu_0-(\rho+\sigma)=\mathbf{v}_{\rho,\sigma}(\{f_0,g_0\})>\mathbf{v}_{\rho,\sigma}(\{f_i,g_j\})=\tau_i+\nu_j-(\rho+\sigma),\] one has $\mathbf{f}_{\rho,\sigma}(\{f,g\})=\{f_0,g_0\}$.

\ep

\noindent \textbf{2.2.} Assume that $f,g\in K[X,Y]$ are $(\rho,\sigma)$-homogeneous of degree $v,u$ respectively, one has
\be\lb{20}\rho X\frac{\pt f}{\pt X}+\sigma Y\frac{\pt f}{\pt Y}=v f\ee and
\be\lb{211}\rho X\frac{\pt g}{\pt X}+\sigma Y\frac{\pt g}{\pt Y}=u g.\ee

 Subtract $\frac{\pt g}{\pt X}$ times \eqref{20} from $\frac{\pt f}{\pt X}$ times \eqref{211}, one has
\be\lb{22}\sigma Y(\frac{\pt f}{\pt X}\frac{\pt g}{\pt Y}-\frac{\pt f}{\pt Y}\frac{\pt g}{\pt X})=ug\frac{\pt f}{\pt X}-vf\frac{\pt g}{\pt X}.\ee
 Subtract $\frac{\pt g}{\pt Y}$ times \eqref{20} from $\frac{\pt f}{\pt Y}$ times \eqref{211},  one has
\be\lb{23} -\rho X(\frac{\pt f}{\pt X}\frac{\pt g}{\pt Y}-\frac{\pt f}{\pt Y}\frac{\pt g}{\pt X})=ug\frac{\pt f}{\pt Y}-vf\frac{\pt g}{\pt Y}.\ee
If $v,u$ are integers (which may be 0), then the above two equations become
\be\lb{12}\sigma Y\{f,g\}=f^{-u+1}g^{v+1}\frac{\pt}{\pt X}(g^{-v}f^u);\ee
\be\lb{13}-\rho X\{f,g\}=f^{-u+1}g^{v+1}\frac{\pt}{\pt Y}(g^{-v}f^u).\ee

The above argument is due to Dixmier; see Lemma 1.4 of \cite{d}.
\ble\lb{80}
Assume that $f,g\in K[X,Y]\setminus K$, and  $(\rho,\sigma)\in \mathbb{Z}^2\setminus\{(0,0)\}$. Assume that $f,g$ are $(\rho,\sigma)$-homogeneous of degree $v,u$ respectively. Then $v,u$ are integers and \[\{f,g\}=0\ \text{if and only if}\ f^u=\lambda g^v\ \text{for some }\ \la\in K^*.\]
\ele
\bp
Assume that $\{f,g\}=0$. Noticing that $f^{-u+1}g^{v+1}\ne 0$, by \eqref{12} and \eqref{13} one has \[\frac{\pt}{\pt X}(g^{-v}f^u)=0\ and\ \frac{\pt}{\pt Y}(g^{-v}f^u)=0,\] which implies  $f^u=\lambda g^v$ for some $\la\in K$. As $f\ne 0$, $\la\in K^*$.

Assume that $f^u=\lambda g^v$ for some $\la\in K^*$. Then by \eqref{12} and \eqref{13} one has \[\sigma Y\{f,g\}=\rho X\{f,g\}=0.\] As at least one of $\sigma$ and $\rho$ is nonzero,  one has $\{f,g\}=0$.
\ep

\begin{prop}\lb{27}
Assume that $f,g\in K[X,Y]\setminus K$, and  $(\rho,\sigma)\in \mathbb{R}^2\setminus\{(0,0)\}$. Assume that $f,g$ are $(\rho,\sigma)$-homogeneous of degree $v,u$ respectively, and $\{f,g\}=0$.

(1) If $\rho,\sigma$ are linearly independent over $\mathbb{Q}$, then $f,g$ are monomials, and $E(f), E(g)$ are on the same ray through O.

(2) If $\rho,\sigma$ are linearly dependent over $\mathbb{Q}$, then one can replace $(\rho,\sigma)$ with some suitable multiple of it so that $(\rho,\sigma)\in \mathbb{Z}^2\setminus\{(0,0)\}$ and  $v,u\in \mathbb{Z}_{\ge 0}$. One has that for some $\lambda\in K^*$, \be \lb{26}f^u=\lambda g^v,\ee and $f,g$ are in one of the following cases:

(2.1) $v=u=0$. In this case $Convex(E(f))$ and $Convex(E(g))$ are on the ray through O orthogonal to $(\rho,\sigma)$;

(2.2) $v,u>0$. Assume that $d=gcd(v,u)$ and ${v}_0=v/d, u_0=u/d $. Then there exists some polynomial $h\in K[X,Y]$ such that $f=\gamma h^{{v}_0}, g=\beta h^{u_0}$, $\gamma\beta\ne 0$. In this case $Convex(E(f))$ and $Convex(E(g))$ are homothetic with respect to O.

In all the 3 cases, one has that $Cone(Convex(E(f)))=Cone(Convex(E(g)))$.

\end{prop}
\bp
(1) Assume that $\rho,\sigma$ are linearly independent over $\mathbb{Q}$.  Then $f,g$ are monomials. It follows from \eqref{54} that $E(f), E(g)$ are on the same ray through O. It is clear that $Cone(Convex(E(f)))=Cone(Convex(E(g)))$.

(2) Assume that  $\rho,\sigma$ are linearly dependent over $\mathbb{Q}$. After replacing  $(\rho,\sigma)$ with $(\tau\rho,\tau\sigma)$ for some suitable $\tau\ne 0$ one has $(\tau\rho,\tau\sigma)\in \mathbb{Z}^2\setminus\{(0,0)\}$, $v,u\in \mathbb{Z}$ and at least one of $v,u$ is $\ge 0$. As $\{f,g\}=0$,  \eqref{12} and \eqref{13} are just
\[ \frac{\pt}{\pt X}(g^{-v}f^u)=0,\ \ \frac{\pt}{\pt Y}(g^{-v}f^u)=0, \] then $g^{-v}f^u=\la\in K$. As $f,g$ are both nonzero polynomials, $\la\ne 0$.

So \[f^u=\lambda g^v, \la\ne 0.\] It is clear that if $u=0$ then $v=0$, and vice versa. So there are 2 cases:

(2.1) $v=u=0$.

In this case $Convex(E(f))$ and $Convex(E(g))$ are both on the line through O orthogonal to $(\rho,\sigma)$.  One has that $Cone(Convex(E(f)))=Cone(Convex(E(g)))$.

(2.2) $v\ne 0, u\ne 0$.

 As $g^{-v}f^{u}$ is a constant, if $u>0>v$ then both $f$ and $g$ are constant polynomials, which contradicts to the assumption. Similarly the case $v>0>u$ can not happen.

So one must have $v>0,u>0$. Let $d=gcd(v,u)$ and $v_0=v/d, u_0=u/d $. As $K[X,Y]$ is a unique factorization domain, one has \[f=\gamma h_1^{s_1}\cdots h_k^{s_k}, g=\beta h_1^{t_1}\cdots h_k^{t_k}, s_i>0,t_i>0,i=1,\cdots, k;\] where $h_i$'s are irreducible monic polynomials. By \eqref{26}, \[ (\gamma h_1^{s_1}\cdots h_k^{s_k})^{u_0 d}= \la (\beta h_1^{t_1}\cdots h_k^{t_k})^{v_0 d}. \] Then $h_i^{s_iu_0d}=h_i^{t_iv_0d}$ and $s_iu_0=t_iv_0,$ for $i=1,\cdots, k$. As $gcd( v_0,u_0)=1$, $v_0$ divides $s_i$. Let $s_i/v_0=m_i, i=1,\cdots, k$. Then $t_i/u_0=m_i, i=1,\cdots, k.$ So \[ s_i=v_0 m_i, t_i=u_0 m_i; i=1,\cdots, k. \] Let $h= h_1^{m_1}\cdots h_k^{m_k}$.
Then $f=\gamma h^{v_0}, g=\beta h^{u_0}$, $\gamma\beta\ne 0$. So $Convex(E(f))$ and $Convex(E(g))$ are homothetic with respect to O, thus $Cone(Convex(E(f)))=Cone(Convex(E(g)))$.

It is verified that in all the 3 cases $Cone(Convex(E(f)))=Cone(Convex(E(g)))$.
\ep

For $f\in K[X,Y]$, let \be\lb{91}C(f)=\{g\in K[X,Y]|\{f,g\}=0\}.\ee  It is easy to verify that $C(f)\supseteq K[f]$.

\ble\lb{79}
  Assume that $f\in K[X,Y]\setminus K$, and that $f$ is $(\rho,\sigma)$-homogeneous of degree $v$ for some $(\rho,\sigma)\in \mathbb{Z}^2\setminus\{(0,0)\}$. Let $g=\sum_{i=0}^k g_i$ be the $(\rho,\sigma)$-homogeneous decomposition of $g\in K[X,Y]\setminus \{0\}$, with $\mathbf{v}_{\rho,\sigma}(g_i)=u_i$ for all $i$. Then $g\in C(f)$ if and only if $g_i\in C(f)$ for all $i$.
\ele
  \bp
By Proposition \ref{78}, $\{f,g_i\}\in V_{\rho,\sigma}(v+u_i-\rho-\sigma)$ for all $i$. One has $\{f,g\}=\sum_{i=0}^k  \{f,g_i\}$. Then $\{f,g\}=0$ if and only if  $\{f,g_i\}=0$ for all $i$.
\ep
For any $f\in K[X,Y]\setminus K$, there exists a unique decomposition \be\lb{73}f=\la h^m,\ee such that

 (1) $\la\ne 0, h\in K[X,Y]$ is monic, and $m\ge 1$;

 (2) $h$ is not a proper power of some polynomial in $K[X,Y]$, i.e. $h$ cannot be written as $h= g^n$ with $ g\in K[X,Y]$ and $n>1$.

  We will refer to such decomposition \eqref{73}  as the \textbf{power decomposition} of $f$.

\bco\lb{90}
Assume that $f\in K[X,Y]\setminus K$, and that $f$ is $(\rho,\sigma)$-homogeneous of degree $v\ne 0$ for some $(\rho,\sigma)\in \mathbb{Z}^2\setminus\{(0,0)\}$. Assume $f=\mu h^m$  is the  power decomposition of $f$. Then

(1) $h$ is also $(\rho,\sigma)$-homogeneous and $C(f)=K[h]$.

(2) $C(f)=K[f]$ if and only if $m=1$.
\eco
\bp
(1)
By $f=\mu h^m$ (with $\mu\ne 0$), $h$ must also be $(\rho,\sigma)$-homogeneous.

Assume that $g\in C(f)\setminus K$ is $(\rho,\sigma)$-homogeneous of degree $u$. Then as we have seen in Proposition \ref{27}, \[f^u=\lambda g^v, for\ some\ \la\ne 0.\] As $f=\mu h^m$, one has \[\mu^u h^{mu}=(\mu h^m)^u=\la g^v.\] By the uniqueness of the power decomposition of $g$, one knows that $g=\gamma h^k$ for some $k\ge 0,\gamma\ne 0$. Thus $g\in K[h]$, and $C(f)\subseteq K[h]$ by Lemma \ref{79}.

As it is clear that $K[h]\subseteq C(f)$, one has $C(f)=K[h]$.

(2) follows from (1).
\ep

\bco\lb{102}
Assume that $f\in K[X,Y]\setminus K$, and that $f$ is $(\rho,\sigma)$-homogeneous of degree $v\ne 0$ for some $(\rho,\sigma)\in \mathbb{Z}^2\setminus\{(0,0)\}$. Assume that $g\in C(f)$. Then

(1) $E(g)\subseteq Cone(Convex(E(f)))$;

(2) If $E(f)\subseteq V_+$ (resp. $E(f)\subseteq V_-$) then $E(g)\subseteq V_+$ (resp. $E(g)\subseteq V_-$).
\eco
\bp
(1) Assume $f=\mu h^m$  ($\mu\ne 0$) is the  power decomposition of $f$. By (1) of Corollary \ref{90}, $h$ is also $(\rho,\sigma)$-homogeneous and $C(f)=K[h]$. As
\[Cone(Convex(E(h)))=Cone(Convex(E(f)))\] and $g\in K[h]$, one has $E(g)\subseteq Cone(Convex(E(f)))$.

 (2) follows from (1).

\ep

\noindent \textbf{2.3.} Recall that $A_1$ is the first Weyl algebra over $K$ generated by $p,q$ with \[[p,q]=1.\] One knows that $A_1$ is a simple Noetherian domain of Gelfand-Kirillov dimension 2. One knows that  $ \{p^i q^j| i,j\in \mathbb{Z}_{\ge 0}\} $ is a basis of $A_1$. Identify $A_1$ with $K[X,Y]$ by \[K[X,Y]\rt A_1, \sum \alpha_{ij}X^iY^j\mapsto \sum\alpha_{ij}p^iq^j.\]
Denote its inverse by  \be\lb{94}\Phi: A_1\rt K[X,Y], p^iq^j\mapsto X^iY^j.\ee We will say that $z\in A_1$ is a monomial (resp. without constant term, etc.) if so is $\Phi(z)$.

For $z\in A_1\setminus\{0\}$ define
\[E(z)=E(\Phi(z))\ \text{and}\  Convex(E(z))=Convex(E(\Phi(z))).\]
For $(\rho,\sigma)\in \mathbb{R}^2\setminus \{(0,0)\}, \rho+\sigma\ge 0$, let \[\mathbf{v}_{\rho,\sigma}(z)=\mathbf{v}_{\rho,\sigma}(\Phi(z))\] and
\[\mathbf{f}_{\rho,\sigma}(z)=\mathbf{f}_{\rho,\sigma}(\Phi(z)).\] The polynomial $\mathbf{f}_{\rho,\sigma}(z)$ is called the polynomial $(\rho,\sigma)$-associated with $z$.

By $[p,q]=1$ one obtains the following result (Lemma 2.1 of \cite{d})

\[q^i p^s=p^s q^i +\sum_{j=1}^{min\{i,s\}}(-1)^j j!
\binom{i}{j}\binom{s}{j} p^{s-j}q^{i-j}.\] Note that the coefficient of $ p^{s-j}q^{i-j}$ is 0 if $j>i$ or $j>s$.
Then it follows that
\be \begin{aligned}p^{s_1}q^{i_1}\cdot p^{s_2}q^{i_2}&=\sum_{j=0}^{min\{i_1,s_2\}}(-1)^j j!
\binom{i_1}{j} \binom{s_2}{j}p^{s_1+s_2-j}q^{i_1+i_2-j}\\
&=p^{s_1+s_2}q^{i_1+i_2}-\binom{i_1}{1}\binom{s_2}{1} p^{s_1+s_2-1}q^{i_1+i_2-1}+\cdots \end{aligned}\ee

Set $min\{i_1,s_2\}=l, min\{i_2,s_1\}=m$. Then

$[p^{s_1}q^{i_1}, p^{s_2}q^{i_2}]$
\bee \begin{aligned}
[p^{s_1}q^{i_1}, p^{s_2}q^{i_2}]
&=\sum_{j=1}^{max{\{l,m\}}}(-1)^j j!\left[\binom{i_1}{j}\binom{s_2}{j}-\binom{i_2}{j}\binom{s_1}{j}\right] p^{s_1+s_2-j}q^{i_1+i_2-j}\\
&=\left |\begin{array}{l}
s_1\  i_1\\
s_2\  i_2
\end{array}\right | p^{s_1+s_2-1}q^{i_1+i_2-1}-(1/2!)\left |\begin{array}{l}
s_1(s_1-1),\  i_1(i_1-1)\\
s_2(s_2-1),\  i_2(i_2-1)
\end{array}\right | p^{s_1+s_2-2}q^{i_1+i_2-2}+\cdots \end{aligned}\eee

The following useful result is due to Dixmier.
\bthm \lb{31}(Lemma 2.7 of \cite{d}, Proposition 3.2 of \cite{jo})\\
Let $ z,w \in A_1\setminus \{0\}$. Assume that $(\rho,\sigma)\in \mathbb{R}^2\setminus \{(0,0)\}$ with $\rho+\sigma> 0$. Set $f=\mathbf{f}_{\rho,\sigma}(z)$ and $g=\mathbf{f}_{\rho,\sigma}(w)$. Then

(1) $\mathbf{f}_{\rho,\sigma}(zw)=fg$;

(2) If $\{f,g\}\ne 0$ then \[\mathbf{f}_{\rho,\sigma}([z,w])=\{f,g\}\ and \ \mathbf{v}_{\rho,\sigma}([z,w])=\mathbf{v}_{\rho,\sigma}(z)+\mathbf{v}_{\rho,\sigma}(w)-(\rho+\sigma);\]

(3) If $\{f,g\}= 0$, then $ \mathbf{v}_{\rho,\sigma}([z,w])<\mathbf{v}_{\rho,\sigma}(z)+\mathbf{v}_{\rho,\sigma}(w)-(\rho+\sigma)$.

\ethm
\bigskip

\noindent \textbf{2.4.} Now we recall the $\mathbb{Z}$-grading on $A_1$.

For any $z\in A_1$, let $ad(z)$ be the derivation \[ad(z):A_1\rt A_1, w\mapsto [z,w]=zw-wz.\]

Then $(ad(pq))(p^i q^s )=[pq, p^i q^s ]=(s-i)p^i q^s$.
Let \[ D_t=Span\{p^i q^s |s-i=t\},\  t\in \mathbb{Z}.\] One knows that $D_0=K[pq]$. If $z\in D_k$, then it is said to be \textbf{homogeneous} of degree $k$. For $i\in \mathbb{Z}$, set
\[D_{\ge i}:=\bigoplus_{i\le t\in \mathbb{Z}} D_t,D_{>i}:=\bigoplus_{i< t\in \mathbb{Z}} D_t.\] $D_{\le i} $ and $D_{<i} $ is defined analogously.

As $D_t\cdot D_s\subseteq D_{t+s}$ for any $t,s\in \mathbb{Z}$, it is clear that $D_{\ge 0}$ and $D_{\le 0}$ are subalgebras of $A_1$.

\ble\lb{32} One has
\[ [D_{\ge 0},D_{\ge 0}]=D_{>0}\ and\ [D_{\le 0},D_{\le 0}]=D_{<0}.\]
\ele
\bp
By $[D_i,D_j]\subseteq D_{i+j}$ for all $i,j\in \mathbb{Z}_{\ge 0}$, and $[D_0,D_0]=0$, one has $[D_{\ge 0},D_{\ge 0}]\subseteq D_{>0}$.

Assume $k>0$. For any $i\ge 0$, as $[pq, p^iq^{i+k}]=kp^iq^{i+k}$, one has $D_k\subseteq [D_{\ge 0},D_{\ge 0}]$. Therefore $[D_{\ge 0},D_{\ge 0}]$ contains all $D_k$ with $k>0$, thus $[D_{\ge 0},D_{\ge 0}]=D_{>0}$.

The equation $[D_{\le 0},D_{\le 0}]=D_{<0}$ is proved similarly.
\ep

It is clear that for any homogeneous element $z$ in $A_1$, $C(z)$ is a graded subalgebra of $A_1$.
\bthm[Theorem 2.2 and Remark 2.6 of \cite{ggv}]\lb{101}

Assume that $z\in D_i\setminus K$, $i\in \mathbb{Z}$.

(1) If $i=0$ then $C(z)=D_0$;  if $i\ne 0$ then $dim(C(z)\cap D_j)\le 1$ for any $j\in \mathbb{Z}$.

(2) If $i> 0$ and $dim(C(z)\cap D_j)= 1$, then $j\ge 0$;  if $i<0$ and \\ $dim(C(z)\cap D_j)= 1$, then $j\le 0$.

\ethm
\bco\lb{51}
Assume that $z\in D_i\setminus\{0\}$. If $i=\pm 1$ then $C(z)=K[z]$.
\eco
\bp We will prove it in the case $i=1$. Note that $C(z)$ and $K[z]$ are both graded subalgebras of $D_{\ge 0}$ with $K[z]\subseteq C(z)$.
As $dim(K[z]\cap D_j)= 1$ for any $j\ge 0$, and $dim(C(z)\cap D_j)\le 1$ for any $j\ge 0$,  one has $K[z]\cap D_j=C(z)\cap D_j$ for any $j\ge 0$, thus $C(z)=K[z]$.
\ep

For any $z\in A_1\setminus\{0\}$, $z=z_0+\cdots+z_m$ with $m\ge 0$ satisfying $z_i\in D_{k_i}\setminus\{0\}$ for $i=0,\cdots,m$ and $k_i>k_{i+1}$ for $i=0,\cdots,m-1$, is called the \textbf{homogeneous decomposition} of $z$ and $z_0$ is called the \textbf{leading component} of $z$. The following result is clear.

\ble\lb{103}
Assume that $z, w\in A_1\setminus\{0\}$ and \[z=z_0+\cdots+z_m, w=w_0+\cdots+w_n  \] are the respective homogeneous decomposition of $z,w$. Then

(1) $z_0w_0$ is the leading component of $zw$.

(2) If $[z_0,w_0]\ne 0$, then $[z_0,w_0]$ is the leading component of $[z,w]$.
\ele

\noindent \textbf{2.5.} Now we describe the automorphism group $\text{Aut}(A_1)$ and $\text{Aut}(K[X,Y])$ of the $K$-algebras $A_1$ and $K[X,Y]$ respectively.

The automorphism group $\text{Aut}(K[X,Y])$ is generated by linear automorphisms
\be \lb{9}X \rt aX+bY,\ \ Y\rt cX+dY\ \ (\left |
\begin{array}{ll}
a & b \\
c & d
\end{array}
\right |\ne 0)\ee and triangular automorphisms
\be\lb{10}X\rt X+f(Y),\ \ Y\rt Y\ \ (f(Y)\in K[Y]).\ee

Recall that for any $\phi\in \text{Aut}(K[X,Y])$, its Jacobian determinant is \[J(\phi)=\left |
\begin{array}{ll}
\partial \phi(X)/\partial X & \partial \phi(X)/\partial Y \\
\partial \phi(Y)/\partial X & \partial \phi(Y)/\partial Y
\end{array}
\right |.\]

The group $\text{Aut}_n (K[X,Y])$ consists of those automorphisms $\phi$ of $K[X,Y]$ with $J(\phi)=1$, which is a normal subgroup of $\text{Aut} (K[X,Y])$. It is generated by the automorphisms in \eqref{9} with $\left |
\begin{array}{ll}
a & b \\
c & d
\end{array}
\right |=1$ and the automorphisms in \eqref{10}.

The polynomial algebra $K[X,Y]$ over $K$ with the Poisson product defined as in \eqref{7} is a Poisson algebra. It is usually called the  1st symplectic Poisson algebra, and is denoted by $S_1(K)$. Note that the $C(f)$ defined in \eqref{91} for $f\in K[X,Y]$ is the centralizer of $f$ in $S_1(K)$.

 Let $\text{Aut}(S_1(K))$ denote the group of automorphisms of $S_1(K)$. An automorphism $\phi$ of the polynomial algebra $K[X,Y]$ is an automorphism of $S_1(K)$ if and only if it satisfies \[\{\phi(f),\phi(g)\}=\{f,g\}, \forall f,g\in K[X,Y]. \] It follows that \[\text{Aut}(S_1(K))\cong \text{Aut}_n(K[X,Y]).\]
\bthm (\cite{d} \cite{ml})
  One has
\[\text{Aut}(A_1)\cong \text{Aut}_n(K[X,Y]).\]
\ethm

If $z_1,z_2\in A_1$ and there exists some $\psi\in \text{Aut}(A_1)$, such that $\psi(z_1)=z_2$, then we say that $z_1$ is conjugate to $z_2$.\\[3mm]

For $\la\in K^*$, let $\psi_\la\in \text{Aut}(A_1)$ be defined such that \[\psi_\la(p)=\la p, \psi_\la(q)=\la^{-1} q.\] Note that $E(\psi_\la(z))=E(z)$ for any $z\in A_1$.

Let $\psi_0\in \text{Aut}(A_1)$ be defined such that  \be\lb{70}\psi_0(p)=q, \psi_0(q)=-p.\ee  Then \[\psi_0^2(p)=-p, \psi_0^2(q)=-q, \] and $\psi_0$ has order 4 in $\text{Aut}(A_1)$.

Under the linear identification \[\Phi:A_1\rt S_1(K), \quad (see \eqref{94}) \]  $\psi_\la$ corresponds to $\tilde{\psi}_\la\in \text{Aut}(S_1(K))$ such that \[\tilde{\psi}_\la(X)=\la X, \tilde{\psi}_\la(Y)=\la^{-1} Y; \] $\psi_0$ corresponds to $\tilde{\psi}_0\in \text{Aut}(S_1(K))$ such that  \be\lb{71}\tilde{\psi}_0(X)=Y, \tilde{\psi}_0(Y)=-X.\ee

Denote the subgroup of $\text{Aut}(A_1)$ consisting of  all the $\psi_\la$ with $\la\in K^*$ by $G_0$. It is clear that $\psi_0$ normalize $G_0$. Call the subgroup of $\text{Aut}(A_1)$ generated by $G_0$ and  $\psi_0$ by $G_1$. Then $G_1$ is the disjoint union of $G_0$ and $G_0\psi_0$.

 Similarly, denote the subgroup of $\text{Aut}(S_1(K))$ consisting all the $\tilde{\psi}_\la$ with $\la\in K^*$ by $\widetilde{G}_0$. One has that $\tilde{\psi}_0$ normalize $\widetilde{G}_0$. Call the subgroup of $\text{Aut}(S_1(K))$  generated by $\widetilde{G}_0$ and  $\tilde{\psi}_0$ by $\widetilde{G}_1$.

 It is clear that there is a unique isomorphism
\be\lb{96}G_1\rt \widetilde{G}_1, \psi\mapsto \tilde{\psi}\ee that maps $\psi_\la\rt \tilde{\psi}_\la, \psi_0\rt \tilde{\psi}_0$, and
\[\tilde{\psi}(\Phi(z))=\Phi(\psi(z)), \forall z\in A_1, \forall \psi\in G_1.  \]

The following result is clear.
\ble One has that \[ \psi_0(D_i)=D_{-i}, i\in \mathbb{Z}.\] And it follows that \[ \psi_0(D_{\ge 0})=D_{\le 0}, \psi_0(D_{\le 0})=D_{\ge 0}.\]\ele

\noindent \textbf{2.6.} Next we define the Newton Polygon for elements in $A_1$.

If  $z\in A_1\setminus \{0\}$, let
\[ NTP(z)=\{(x,y)\in (\mathbb{R}_{\ge 0})^2|\exists  t\in \mathbb{R}_{\ge 0}, (x,y)+t(1,1)\in Convex(E(z)) \},\]
which is called the (solid) \textbf{Newton polygon} of $z$. If $z=0$, then $E(z)$ and $NTP(z)$ are both defined to be the empty set.

Let $(\rho,\sigma)$ runs continuously from $(1,-1)$ to $(-1,1)$ on the path defined by
\[ (\rho(t), \sigma(t))=\left\{
\begin{array}{ll}
(1, t),    & t\in [-1,1]; \\
(2-t, 1),  & t\in [1,3].
\end{array}
\right. \]
The set
\[Roof(z)=\bigcup_{-1<t<3} Convex(E(\mathbf{f}_{\rho(t),\sigma(t)}(z))) \] is called the roof of $NTP(z)$, which is always concave. Note that the set $Convex(E(\mathbf{f}_{\rho(t),\sigma(t)}(z)))$ is either a vertex or an edge of $NTP(z)$. If $|E(\mathbf{f}_{\rho(t),\sigma(t)}(z))|=1$, then $Convex(E(\mathbf{f}_{\rho(t),\sigma(t)}(z)))$ is a vertex. If $|E(\mathbf{f}_{\rho(t),\sigma(t)}(z))|\ge 2$, then $Convex(E(\mathbf{f}_{\rho,\sigma}(z)))$ is an edge.

Note that \[ NTP(z)=\{(x,y)\in (\mathbb{R}_{\ge 0})^2|\exists  t\in \mathbb{R}_{\ge 0}, (x,y)+t(1,1)\in Roof(z)\},\] so $NTP(z)$ is determined by $Roof(z)$.

Let $z=\al_1 p+\al_2 p^2q^3+\al_3 p^3q+\al_4 p^4q^2+\al_5p^5$ with $\al_i\ne 0, i=1,\cdots, 5$. Then $NTP(z)$ is the pentagon $OP_1P_2P_3P_4$ and $Roof(z)$ is the polygonal chain: $(P_1,P_2,P_3)$. See Figure 1.

\begin{figure} [htp]
  \centering

    \begin{tikzpicture}[thick,>=stealth]
      \draw [->] (0,0) -- (6,0) node [below] {$x$};
      \draw [->] (0,0) -- (0,4) node [left] {$y$};

      \coordinate[label=above:$P_1$] (p1) at (5,0);
      \coordinate[label=above right:$P_2$] (p2) at (4,2);
      \coordinate[label=above:$P_3$] (p3) at (2,3);
      \coordinate[label=left:$P_4$] (p4) at (0,1);

      \draw (p1) -- (p2) -- (p3) -- (p4);
      \foreach \j in {1,2,3,4}
      {
       \fill (p\j) circle (1.5pt);
      }
      \foreach \p in {(0,0),(1,0),(3,1)}
      {
        \fill \p circle (1.5pt);
      }
      \node [left] at (3,2) {$z$};
      %\node [left] at (0,3) {$3$};
      \node [below] at (1,0) {$1$};
      \node [below] at (5,0) {$5$};
      \node [below] at (0,0) {$O$};
    \end{tikzpicture}
   \caption*{Figure 1}
\end{figure}

Let $a=min\{j-i|(i,j)\in E(f) \}, b=max\{j-i|(i,j)\in E(f) \}$. Then $Roof(z)$ and $NTP(z)$ are both in the region $\{(x,y)\in (\mathbb{R}_{\ge 0})^2|a\le y-x\le b \}.$

There are 2 cases of $Roof(z)$ with $z\ne 0$.

(1) $\mathbf{f}_{\rho(t), \sigma(t)}(z)$ is a monomial, for any $t\in (-1,3).$ Now $z$ is in some single homogeneous space $D_k$. One has \[Roof(z)=E(\mathbf{f}_{1,1}(z)),\] which is just a point. Otherwise, it is in next case.

(2) There exists some $n\ge 1$ and $t_1,t_2,\cdots,t_n$, where $-1=t_0<t_1<t_2<\cdots<t_n<t_{n+1}=3$, are all the $t\in (-1,3)$ such that $\mathbf{f}_{\rho(t_i), \sigma(t_i)}(z)$ is not a monomial. One has \[Roof(z)=\bigcup_{i=1}^n Convex(E(\mathbf{f}_{\rho(t_i),\sigma(t_i)}(z))), \] which consists of $n$ edges.

 The following observation is clear and we omit its proof.
\ble\lb{60}
The following statements are equivalent:

(1) $z\in D_{\ge 0}$;  (2) $E(z)\subseteq V^+$; (3) $Roof(z)\subseteq V^+$;  (4) $Cone(Roof(z))\subseteq V_+$; \\(5) $NTP(z)\subseteq V_+$.
\ele An analogous result for $z\in D_{\le 0}$ also holds.\\[3mm]

\section{Proof of the main result}
 \setcounter{equation}{0}\setcounter{theorem}{0}

Let \[{\Gamma}=\{(z,w)\in A_1^2|[z,w]=1\}.\]
For any $\psi\in \text{Aut}(A_1)$ and $(z,w)\in \Gamma$, $[\psi(z), \psi(w)]=1$ thus $(\psi(z), \psi(w))\in {\Gamma}.$ So $\text{Aut}(A_1)$ acts on ${\Gamma}$ (faithfully) by
\[ \psi\cdot (z,w)= (\psi(z), \psi(w)).\]

If $(z,w)\in \Gamma$ then $(w,-z)\in \Gamma$. Let \[ \eta:\Gamma\rt \Gamma, (z,w)\mapsto (w,-z). \] Let $U$ be the group of transformations of $\Gamma$ generated by $G_1$ (defined in \textbf{2.5}) and $\eta$. It is directly verified that
\[ {\eta}{\psi}= {\psi}{\eta}, \forall {\psi}\in G_1.\]

Let $\Omega$ be the collection of $(f,g)$ in $K[X,Y]^2$ such that

(1) $\{f,g\}=1$;

(2) $f,g$ are $(\rho,\sigma)$-homogeneous for some $(\rho,\sigma)\in \mathbb{R}^2\setminus \{(0,0)\}$.

For $(f,g)\in \Omega$ and $\tilde{\psi}\in \widetilde{G}_1$ (defined in \textbf{2.5}), one has  $(\tilde{\psi}(f), \tilde{\psi}(g))\in \Omega.$ So $\widetilde{G}_1$ acts on $\Omega$ (faithfully) by
\[ \tilde{\psi}\cdot (f,g)= (\tilde{\psi}(f), \tilde{\psi}(g)).\]

If  $(f,g)\in \Omega$ then $(g,-f)\in \Omega$.  Let \[\tilde{\eta}:\Omega\rt \Omega, (f,g)\mapsto (g,-f). \] Let $\widetilde{U}$ be the group of transformations of $ \Omega$ generated by $\widetilde{G}_1$ and $\tilde{\eta}$. It is directly verified that
\[ \tilde{\eta}\tilde{\psi}= \tilde{\psi}\tilde{\eta}, \forall \tilde{\psi}\in \widetilde{G}_1.\]

\begin{prop}\lb{30}

Assume $(f,g)\in \Omega$. Then up to the action of $\widetilde{U}$, $(f,g)$ will be in one and only one of the following cases:

(1) $(X,Y)$;

(2) $(\alpha X+\beta Y, \gamma X+\delta Y), \alpha,\delta,\beta,\gamma\in K, \alpha\delta-\beta\gamma=1, \alpha\delta\beta\gamma\ne 0$;

(3) $(X+\lambda Y^n, Y), \la\ne 0, n\ge 1$;

%(4) $(X, Y+\mu X^n), \mu\ne 0, n\ge 1$.

(4)  $(X+\la,Y), \la\ne 0$.
\end{prop}
Let us introduce a terminology before the proof. Assume that $f,h\in K[X,Y]$, $f=\sum \alpha_{ij} X^iY^j$,  $h=\sum \beta_{ij} X^iY^j$. If whenever $\beta_{ij}\ne 0$, one has $\alpha_{ij}=\beta_{ij}$, then we say that $f$ \textbf{contains} $h$.

\bp
Assume that $(f,g)\in \Omega$. Then $\{f,g\}=1$ and $f,g$ are  $(\rho,\sigma)$-homogeneous for some $(\rho,\sigma)\in \mathbb{R}^2\setminus \{(0,0)\}$ with $\rho+\sigma\ge 0$.

There are 2 cases: (i) Both $f$ and $g$ are monomials; (ii) $f$ or $g$ is not a monomial. In Case (ii), if the 2nd polynomial $g$ of $(f,g)$ is not a monomial,  then after applying $\tilde{\eta}$ one can always assume that the 1st polynomial $f$  of $(f,g)$ is not a monomial.

(i) Suppose that $f, g$ are both monomials.

 Assume that $f=\lambda X^iY^j,g=\mu X^kY^l;$ $\la,\mu\in K^*$. By \eqref{54}, $
\{X^iY^j,X^k Y^l\}=\left |
\begin{array}{ll}
i & j \\
k & l
\end{array}
\right | X^{i+k-1}Y^{j+l-1}.
$  Then \[i+k=1, j+l=1, il-kj\ne 0.\] Thus $(i,j)=(1,0), (k,l)=(0,1)$ or vice verse. One has \[(f,g)=(\lambda X,\lambda^{-1}Y), or\ (f,g)=(\lambda Y,-\lambda^{-1}X); \lambda\in K^*.\] After applying some transformation in $\widetilde{G}_1$, $(f,g)=(X,Y)$, which is in Case (1).

(ii)  Assume that $f$ is not a monomial.

Now $\rho,\sigma$ are linearly dependent over $\mathbb{Q}$. One can assume that $(\rho,\sigma)=(r,s)$, where $r,s$ are coprime integers and $r+s\ge 0$.  Assume that  $f,g$ are  $(r,s)$-homogeneous of degree $t,u$ respectively. By $\{f,g\}=1$ one has $t+u-(r+s)=0$ thus \be\lb{47} t+u=r+s=\mathbf{v}_{r,s}(XY).\ee There are 3 cases to be dealt with:

(a) $r, s>0$; (b) $(r,s)=(0,1)$; (c) $r>0>s$.

The remaining cases $(r,s)=(1,0)$ and $s>0>r$ can be transformed to Case (b) and (c) respectively after applying $\tilde{\psi}_0$ (defined in \eqref{71}).\\[1mm]

(a) $r, s>0$.

If $t=0$ (resp. $u=0$), then $f$ (resp. $g$) is a constant, which contradicts to $\{f,g\}=1$. So $t,u>0$. By \eqref{47}, one has $t=\mathbf{v}_{r,s}(f)<\mathbf{v}_{r,s}(XY)$.

If $\mathbf{v}_{r,s}(X^iY^j)<\mathbf{v}_{r,s}(XY)$, then $ri+sj<r+s$, thus $i=0$ or $j=0$. So \[f=\lambda X+\mu Y^n, \ or,\ f=\beta Y+\gamma X^n; \lambda,\mu,\beta,\gamma\in K^*, n\ge 1.\] See Figure 2.

\begin{figure}[htp]
  \centering
  \begin{tikzpicture}[scale=0.8]
    \tikzstyle{every node}=[font=\small,scale=0.8]
    \draw [->] (0,0) -- (5,0) node [below] {$x$};
    \draw [->] (0,0) -- (0,5) node [left] {$y$};
    %  \node [below left] at (0,0) {$O$};

    \coordinate (w1) at (0,3);
    \coordinate (w2) at (1,0);
    \coordinate (z1) at (0,0);
    \coordinate (z2) at (1,1);
    \coordinate (x1) at (4,4);
    \coordinate (x2) at (1,3);

    \draw (w1) -- node [left] {$f$} (w2);

    \draw [<-] (x1) -- node [below right] {$\left(r,s\right) $} (x2);
    \foreach \j in {1,2}
    {
      \foreach \k in {w,z}
      {
        \fill (\k\j) circle (1.5pt);
      }
    }
    \fill (0,1) circle (1.5pt);
    \fill (x2) circle (1.5pt);
    \node [left] at (0,3) {$3$};
    \node [below] at (1,0) {$1$};
    \node [left] at (0,1) {$1$};
    \node [below] at (0,0) {$O$};
  \end{tikzpicture}
  \caption*{Figure 2}
\end{figure}

 After applying $\tilde{\psi}_0$ if needed, one can assume that $f=\lambda X+\mu Y^n, \la\mu\ne 0$. Then $(r,s)=(n,1)$ and $t=\mathbf{v}_{r,s}(f)=n$, thus $u=r+s-n=1$.

Assume that $g=\sum \gamma_{ij} X^iY^j$. Then \[\mathbf{v}_{r,s}(X^iY^j)=ri+sj=ni+j=1\] for those $i,j$ with $\gamma_{ij}\ne 0$.

If $n>1$ then $i=0,j=1$ and $g=\gamma Y, \gamma\ne 0$.\[1=\{f,g\}=\{\lambda X+\mu Y^n,\gamma Y\}=\lambda \gamma. \]After applying some $\psi\in \widetilde{G}_1$, one can assume that $(f,g)=(X+\beta Y^n,Y), \beta\ne 0$. This is in Case (3).

If $n=1$ then $i=0,j=1$ or $i=1,j=0$. Now $g=\gamma X+\delta Y$. By $\{f,g\}=1$, one has $\lambda\delta-\mu\gamma=1$.  It is in Case (2) if $\delta\gamma\ne 0$; it is in Case (3) if one of  $\delta$ and $\gamma$ is 0.\\[1mm]

(b) $(r,s)=(0,1)$. Now $t+u=r+s=1, t,u\ge 0$. Then $t=1,u=0$ or $t=0,u=1$.

 (b.1) $t=1,u=0$. Assume that $f=F(X)Y, g=H(X)$, where $F,H\in K[X]$. By $\{f,g\}=1$, one has
\[ 1=\{F(X)Y, H(X)\}=  -F(X)H'(X).  \] Then $F(X)$ is a nonzero constant $\al$ and $(f,g)=(\al Y,-\al^{-1} X+r)$. Applying $\tilde{\eta}$, $(f,g)$ becomes $(-\al^{-1} X+r,-\al Y )$. Then applying some transformation in $\widetilde{G}_0$, $(f,g)$ becomes $(X+r,Y )$, which is in Case (4).

 (b.2) $t=0,u=1$. Then $(g,-f)=\tilde{\eta}(f,g)$ is in the situation of (b.1).\\[1mm]

(c) $r>0>s$. As $\{f,g\}=1$,  $f$ and $g$ contain the monomials $\alpha_0 X$ and $\beta_0 Y$ respectively, where $\alpha_0\beta_0\ne 0$. After applying $\tilde{\psi}_0$ if necessary, we assume that  $f$ contains $\alpha_0 X$ and that $g$ contains $\beta_0 Y$. Then, as $r,s$ are coprime, one has

\[f=X\sum_{i=0}^k \alpha_i (X^{-s}Y^r)^i, k\ge 0,  \alpha_k\ne 0;\] and
\[g=Y\sum_{i=0}^l \beta_i (X^{-s}Y^r)^i, l\ge 0, \beta_l\ne 0.\] See Figure 3.

\begin{figure}[htb]
  \centering
\begin{tikzpicture}[scale=0.7]
  \tikzstyle{every node}=[font=\small,scale=0.7]
  \draw [->] (0,0) -- (6,0) node [below] {$x$};
  \draw [->] (0,0) -- (0,6) node [left] {$y$};
  %  \node [below left] at (0,0) {$O$};

  \coordinate (w1) at (1,0);
  \coordinate (w2) at (3,1);
  \coordinate (z1) at (0,1);
  \coordinate (z2) at (2,2);
  \coordinate (x1) at (2,3);
  \coordinate (x2) at (1,5);

  \draw (w1) -- node [below] {$g$} (w2)
  (z1) -- node [above] {$f$} (z2);
  \draw [<-] (x1) -- node [above right] {$ \left(r,s\right)$} (x2);
  \foreach \j in {1,2}
  {
    \foreach \k in {w,z}
    {
      \fill (\k\j) circle (1.5pt);
    }
  }
  \fill (0,0) circle (1.5pt);
  \fill (x2) circle (1.5pt);

  \node [left] at (0,1) {$1$};
  \node [below] at (1,0) {$1$};
  \node [below] at (0,0) {$O$};
\end{tikzpicture}
\caption*{Figure 3}
\end{figure}

One has
\[\mathbf{f}_{-s,r}(f)=  \alpha_k X^{-sk+1}Y^{rk}, \mathbf{f}_{-s,r}(g)=  \beta_l X^{-sl}Y^{rl+1}.\]

By $r>0>s$, it is clear that \[\{\mathbf{f}_{-s,r}(f), \mathbf{f}_{-s,r}(g)\}=(rl-sk+1)\al_k\beta_l X^{-s(k+l)}Y^{r(k+l)}\ne 0,\] thus by Proposition\ \ref{24}, \[\{\mathbf{f}_{-s,r}(f), \mathbf{f}_{-s,r}(g)\}=\mathbf{f}_{-s,r}(\{f,g\})=1,\] which implies $k=l=0$. So $(f,g)$ is still in Case (1).

 Note that for any $\psi\in \widetilde{G}_1$ and $h\in K[X,Y]$, $|E(\psi(h))|=|E(h)|$.
We have shown that $(f,g)$ will be in one of the 4 cases up to the action of $\widetilde{U}$. Since $(|E(f)|,|E(g)|)$ is respectively $(1,1), (2,2),(2,1), (2,1)$ in the 4 cases, we only need to show that  the transformations in $\widetilde{U}$ will not take the $(f,g)$  in Case (3) into  Case (4) but this is clear as there does not exist $\tilde{\psi}\in  \widetilde{G}_1$ with $\tilde{\psi}(X+\la Y^n)=X+\al$, where $n\ge 1$.
\ep
The following observation follows easily from the above proposition.
\bco\lb{100}
Assume $(f,g)\in \Omega$. Then $f$ contains the term $\al X$ or $\al Y$ for some $\al\ne 0$.  If $f$ contains the term $\al X$, then either $f=\al X$, or $f=\al X+\la Y^n$ with $\la\ne 0$ and $n\ge 0$.
\eco

\begin{prop} \lb{3}
Assume that $z,w\in A_1$ and $[z,w]=1$. If $\mathbf{v}_{0,1}(z)\le 1$, then $z$ and $w$ generate $A_1$. Specifically,

(1) If $\mathbf{v}_{0,1}(z)=1$, i.e., $z=f(p)q+g(p)$ for some $f,g\in K[X]$ and $f\ne 0$, then
\be\lb{45}z=\alpha q+g(p), w=\gamma-\alpha^{-1}p+h(z); \alpha\ne 0, h\in K[X].\ee

(2) If $\mathbf{v}_{0,1}(z)=0$, i.e. $z=f(p)$  for some $f\in K[X]$, then \be\lb{46}z=\alpha p+\beta, w={\alpha}^{-1}q+g(p), \alpha\ne 0, g\in K[X].\ee

\end{prop}
\bp
(1) Assume that $\mathbf{v}_{0,1}(w)=j$ and $\mathbf{f}_{0,1}(w)=g(X)Y^j$. Note that $\mathbf{f}_{0,1}(z)=f(X)Y$.

 If $\{f(X)Y,g(X)Y^j\}\ne 0$ then \[j=\mathbf{v}_{0,1}([z,w])=\mathbf{v}_{0,1}(1)=0.\]

 So, if $j>0$ then $\{f(X)Y,g(X)Y^j\}=0$. By Lemma \ref{80}, \[g(X)Y^j= \mu (f(X)Y)^j=\mu f(X)^j Y^j,\  for\ some\ \mu\in K^*.\] Then $g(X)=\mu f(X)^j$. Let $w_1=w-\mu z^j$. One has $[z,w_1]=1$ and $0\le \mathbf{v}_{0,1}(w_1)<j$. If  $\mathbf{v}_{0,1}(w_1)>0$ then repeat the above process until we get some polynomial $h(X)$ such that the element $w'=w-h(z)$ satisfies $\mathbf{v}_{0,1}(w')=0$ and $[z,w']=1$.

 If $j=0$ then let $w'=w$.

Assume $w'=l(p)$ for some $l\in K[X]$. As \[[z,w']=[f(p)q+g(p), l(p)]=f(p)[q, l(p)]=-f(p)l'(p)=1,\] one has $f(p)=\alpha\ne 0, l'(p)=-\al^{-1}$. So $z=\alpha q+g(p), w'=l(p)=-\alpha^{-1} p+\gamma$, from which one gets \eqref{45}. It is clear that $z,w$ generate $A_1$.

(2) Assume that $\mathbf{v}_{0,1}(w)=j$. The assumption that $[z,w]=1$ implies that $j>0$. One has \[0=\mathbf{v}_{0,1}(1)=\mathbf{v}_{0,1}([z,w])=j-1.\] So $j=1$. Assume that $ w=l(p)q+h(p)$ for some $l,h\in K[X]$. Then by $[z,w]=1$ and similar computation as in (1), one gets \eqref{46}. It is clear that $z,w$ generate $A_1$.

\ep

\ble \lb{1}
Assume that $z\in D_{\le 0}$ or $z\in D_{\ge 0}$, and $(i,i)$ with $i\ge 1$ is a vertex of $NTP(z)$, i.e., $\mathbf{f}_{\rho,\sigma}(z)=\lambda X^iY^i$ for some $(\rho,\sigma)$ with $\rho+\sigma>0$ and $\la\ne 0$, then for any $w\in A_1$, $[z,w]\ne 1$.
\ele
\bp We will prove it in the case $z\in D_{\le 0}$. Then there exists some $\epsilon\in \mathbb{R}$ with $0<\epsilon<1$ such that $\mathbf{f}_{-1+\epsilon,1}(z)=\lambda X^iY^i$. Suppose that for some $w\in A_1$, $[z,w]=1$. Assume that $\mathbf{f}_{-1+\epsilon,1}(w)=g(X,Y)$.

 If $w\nsubseteq D_{\le 0}$, then $E(g)\cap \{(x,y)\in \mathbb{R}^2|y>x\ge 0\}\ne \varnothing$ and \[\{\mathbf{f}_{-1+\epsilon,1}(z),\mathbf{f}_{-1+\epsilon,1}(w)\}=\{\lambda X^iY^i, g(X,Y)\},\] which cannot be 1 by Proposition  \ref{30}, and cannot be 0 by Proposition  \ref{27}, thus $[z,w]\ne 1$.

  If $w\subseteq D_{\le 0}$ then $[z,w]\subseteq D_{<0}$ by Lemma \ref{32}, thus $[z,w]\ne 1$.
\ep
Let $\widetilde{A_1}$ be the $K$-algebra generated by $p,q,q^{-1}$ subject to the relation \[q q^{-1}=q^{-1} q=1, [p,q]=1.  \] It is easy to verify that $\{p^iq^k| i\in \mathbb{Z}_{\ge 0}, k\in \mathbb{Z}\}$ is a basis of $\widetilde{A_1}$. Then $A_1$ is the subalgebra of $\widetilde{A_1}$ generated by $p$ and $q$.

It is easy to verify that
 \be\lb{99} q^k \cdot pq=(pq-k) q^k\ee for any $k\in \mathbb{Z}$ as in Section 3.2 of \cite{d}, and it follows that
\be\lb{40}q^k f(pq)=f(pq-k)q^k\ee for any polynomial $f\in K[X]$ and any $k\in \mathbb{Z}$.
\ble\lb{2}
Assume that $z,w\in A_1$ and $[z,w]=1$. Then

If $z$ is contained in some homogeneous space $D_k$ with $k\in \mathbb{Z}$, then \be \lb{42}z=\lambda q, w=\mu p+l(q), \lambda\mu=-1,l\in K[X];\ee or
\be\lb{43}z=\lambda p, w=\mu q+l(p), \lambda\mu=1,l\in K[X].\ee
\ele
\bp
By Lemma \ref{1}, if $z$ is contained in $D_0$ then $[z,w]\ne 1$ for any $w\in A_1$.

Assume that $z\in D_k$ with $k>0$. Then there exists $f,g\in K[X]\setminus \{0\}$ such that \[z=f(pq)q^k, and\ w=w_0 +w_1, w_0=g(pq)q^{-k}, w_1\in \bigoplus_{i\ne -k} D_i.\] Then \be\lb{41}[z,w_0]=1, and\ \ [z,w_1]=0.\ee
Let $h(X)=f(X)g(X-k)$. By \eqref{99} one has

\bee\begin{split} 1&=[f(pq)q^k , g(pq)q^{-k} ]\\&=f(pq) (q^k g(pq)) q^{-k} - g(pq)(q^{-k}f(pq))q^k\\& = f(pq)g(pq-k) q^k\cdot q^{-k}-g(pq)f(pq+k) q^{-k}\cdot q^k\ \ (by\ \eqref{40})\\&=f(pq)g(pq-k)-f(pq+k)g(pq)\\&=h(pq)-h(pq+k).\end{split}\eee Then $h(X)-h(X+k)=1$. So $deg(h(X))=deg(f(X))+deg(g(X))=1$.

As $k>0$, if $deg(g(X))=0$ then $w_0=\lambda q^{-k}$ (for some $\la\ne 0$) is not in $A_1$. So $deg(g(X))=1$ and $deg(f(X))=0$. Then $z=\lambda q^k$ and $w_0=(\mu\ pq+\gamma)q^{-k}$ with $\lambda,\mu\ne 0$. As $w_0\in A_1$, $k=1$ and $\gamma=0$. Then $z=\lambda q, w_0=\mu p, \lambda\mu=-1$. By \eqref{41}, $w_1\in C(\lambda q)=K[q]$. This is the case \eqref{42}.

If $z\in D_k$ with $k<0$, then it will be in the case \eqref{43}.
\ep

\ble\lb{44}
Assume that

(1) $z , w\in A_1\setminus K$, and $z,w$ are both without constant term;

(2) for any $(\rho,\sigma)\in \mathbb{R}^2$ with $\rho+\sigma> 0$, $\{\mathbf{f}_{\rho,\sigma}(z),\mathbf{f}_{\rho,\sigma}(w)\}=0$.

Then \[Cone(Roof(z))=Cone(Roof(w)).\]
\ele
\bp
One has

 \[ Roof(z)=\bigcup_{\tiny{\begin{array}{c} (\rho,\sigma)\\ \rho+\sigma>0\end{array}}}\ Convex(E(\mathbf{f}_{\rho,\sigma}(z))).\]

By the assumption $\{\mathbf{f}_{\rho,\sigma}(z),\mathbf{f}_{\rho,\sigma}(w)\}=0$ for any $(\rho,\sigma)$ with $\rho+\sigma> 0$, applying
Proposition \ref{27} one has, \bee\begin{aligned} Cone(Roof(z))&=\bigcup_{\tiny{\begin{array}{c} (\rho,\sigma)\\ \rho+\sigma>0\end{array}}}\ Cone(Convex(E(\mathbf{f}_{\rho,\sigma}(z))))\\
&=\bigcup_{\tiny{\begin{array}{c} (\rho,\sigma)\\ \rho+\sigma>0\end{array}}}\ Cone(Convex(E(\mathbf{f}_{\rho,\sigma}(w))))\\
&=Cone(Roof(w)).
\end{aligned}\eee
\ep

We are ready to prove the following main result.

\bthm\lb{33}
Assume that $z,w\in A_1,$ $[z,w]=1$, and $z$  is conjugate to some element in $D_{\ge 0}$ (or in $D_{\le 0}$). Then $z$ and $w$ generate $A_1$.
\ethm
\bp
As any element in $D_{\ge 0}$ is conjugate to some element in $D_{\le 0}$ by $\psi_0$ in \eqref{70}, without loss of generality, we assume that $z\in D_{\le 0}$.

If $z',w'\in A_1$, $z'-z\in K,w'-w\in K$, then it is clear that $[z',w']=1$ and $z, w$ generate $A_1$ if and only if $z', w'$ generate $A_1$. So we assume that---
\textbf{both $z$ and $w$ have constant term 0.}\\[2mm]

(i) If $z$ is a monomial, then it follows from Lemma \ref{2} that
\[ z=\lambda p, w=\lambda^{-1} q+f(p), \ \la\ne 0, \ \ f\in K[X].\]
Then $z$ and $w$ generate $A_1$.  If $w$ is a monomial, then one can prove that $z$ and $w$ generate $A_1$ similarly.

(ii) Next we deal with the case that $z$ is not a monomial. (We will show that $[z,w]\ne 1$.) The case that $w$ is not a monomial can be treated similarly.

As $z\in D_{\le 0}$, $Roof(z)\subseteq V_-, NTP(z)\subseteq V_-$. As by assumption $[z,w]=1$,  for any $(\rho,\sigma) $ with $\rho+\sigma> 0$, $\{\mathbf{f}_{\rho,\sigma}(z), \mathbf{f}_{\rho,\sigma}(w)\}=0$ or 1, by Theorem \ref{31}.

If for any $(\rho,\sigma) $ with $\rho+\sigma> 0$, $\{\mathbf{f}_{\rho,\sigma}(z), \mathbf{f}_{\rho,\sigma}(w)\}=0$, then by Lemma \ref{44}, $Cone(Roof(w))=Cone(Roof(z))\subseteq V_-$, thus $Roof(w)\subseteq V_-$. By Lemma \ref{60},  $NTP(w) \subseteq V_-$, thus $w\in A_1^{-}$. Then by Lemma \ref{32}, $[z,w]\in D_{<0}$ and $[z,w]\ne 1$.

In the rest of the proof we assume that for some $(\rho,\sigma) $ with $\rho+\sigma> 0$, $\{\mathbf{f}_{\rho,\sigma}(z), \mathbf{f}_{\rho,\sigma}(w)\}=1$. Set $f=\mathbf{f}_{\rho,\sigma}(z), g= \mathbf{f}_{\rho,\sigma}(w) $. As $z\in D_{\le 0}$, $f$ contains some $\al X$ with $\al\ne 0$, and any $(i,i)$ with $i\ge 1$ is not in $E(z)$ by Lemma \ref{1}. Since $z$ has constant term 0, \[E(z)\subseteq \{(x,y)\in\mathbb{R}^2| x-y\ge 1, y\ge 0  \}, (1,0)\in E(z).\]

As $z$ is not a monomial and $z\in D_{<0}$, there are 2 cases:

(a.1) There exists some unique $(r,s) $ up to a positive constant such that $r\le 0, r+s> 0$, $\mathbf{f}_{r,s}(z)$ contains $\alpha X$, and $\mathbf{f}_{r,s}(z)$ is not a monomial. See Figure 4.

\begin{figure}[htp]
  \centering
  \begin{tikzpicture}[scale=0.8]
    \tikzstyle{every node}=[font=\small,scale=0.8]
    \draw [->] (0,0) -- (5,0) node [below] {$x$};
    \draw [->] (0,0) -- (0,5) node [left] {$y$};

    \coordinate[label=above right:$w$] (w1) at (0,1);
    \coordinate (w2) at (4,0);
    \coordinate (z1) at (4,2);
    \coordinate (z2) at (1,0);
    \coordinate (x1) at (0,0);
    \coordinate (x2) at (4,4);

    \draw (w1) -- (w2) (z1) -- (z2);
    \draw [dashed] (x1) -- (x2);
    \foreach \j in {1,2}
    {
      \foreach \k in {w,z}
      {
        \fill (\k\j) circle (1.5pt);
      }
    }

  \fill (x1) circle (1.5pt);
    \node [below ] at (0,0) {$O$};
    \node [above left=0.5em] at (4,1) {$z$};
    \node [left] at (0,1) {$1$};
    \node [below] at (1,0) {$1$};
  \end{tikzpicture}
  \caption*{Figure 4}
  \end{figure}

(a.2) $\mathbf{f}_{-1,1}(z)$ contains $\alpha X$ and is not a monomial.

We first deal with (a.1). Let $\mathbf{f}_{r,s}(z)=f$ and $\mathbf{f}_{r,s}(w)=g$. By the assumption $[z,w]=1$, one has $E(g)\nsubseteq V_-$. While $E(f)\subseteq V_-$, $\{f, g\}\ne 0$ by Corollary \ref{102}. One has $\{f, g\}\ne 1$ by Corollary \ref{100}. So $[z,w]\ne 1$ and this case cannot happen.

Then we deal with (a.2). Write $z=z_0+z_1, w=w_0+w_1$, where $z_0$, $w_0$ are the respective leading component of $z,w$.

Assume that $w_0\in D_i$. As by assumption $[z,w]=1$, one has $i>0$. Then $[z_0,w_0]\ne 0$ by Theorem \ref{101}, thus $[z_0,w_0]$ is the leading component of $[z,w]$ and $[z_0,w_0]=1$. By \eqref{43} of Lemma \ref{2}, one has $z_0=\al p$ with $\al\ne 0$, which contradicts to the assumption that $\mathbf{f}_{-1,1}(z)$ (which equals $\mathbf{f}_{-1,1}(z_0)=\al X$) is not a monomial.

 Thus $[z,w]\ne 1$.

So we have shown that if $z\in D_{\le 0}$ and $[z,w]=1$ then either $z$ or $w$ is a monomial, which is in the situation of (i) and $z,w$ generate $A_1$.

\ep
\begin{rem}\lb{75}
Assume that  $z\in D_{\le 0}$. If there exists some $w\in A_1$ with $[z,w]=1$, then in the above proof one knows that
\[ z=\lambda p+\gamma, w=\lambda^{-1} q+f(p), \la\ne 0, \ \ f\in K[X].\]

\end{rem}
The following result gives an equivalent formulation of DC.
\begin{prop}
DC holds if and only if for any $z,w\in A_1$ with $[z,w]=1$, there exists some $\psi\in \text{Aut}(A_1)$ such that $\psi(z)\in D_{\le 0}$.
\end{prop}
\bp
The 'if' part is obvious. So we only need to prove the 'only if' part. Assume that DC holds. Let $z,w\in A_1$ such that $[z,w]=1$. Then $z$ and $w$ generate $A_1$. There exists a unique $K$-algebra homomorphism $\psi:A_1\rt A_1$ with $\psi(z)=p,\psi(w)=q$. It is clear that $\psi$ is surjective. The kernel of $\psi$ is an ideal of the simple algebra $A_1$, thus must be the zero ideal. So $\psi\in \text{Aut}(A_1)$ and $\psi(z)\in D_{\le 0}$.

\ep

If $z,w\in A_1$ with $[z,w]=1$, then $z,w$ are both nilpotent in the sense of \cite{d}. It is known that a nilpotent element may not be conjugate to some element in $D_{\ge 0}$ by an automorphism in $A_1$ even when $K$ is algebraically closed; see the remark after Theorem 4.2 of \cite{jo}. But it is not known if any $z\in A_1$ satisfying $[z,w]=1$ for some $w\in A_1$ is conjugate to some element in $D_{\ge 0}$.

\bco
Assume that $z,w\in A_1$ with $[z,w]=1$.
 If $z\in D_{\ge -s}, s>0$, $z=z_{-s}+z'$ with $z_{-s}\in D_{-s}\setminus\{0\}$, $z'\in D_{> -s}$, and $C(z_{-s})=K[z_{-s}]$,
then $z$ and $w$ generate $A_1$.
\eco
\bp
If $w\in D_{\ge 0}$, then the result holds by Theorem \ref{33}. So assume that $w\notin D_{\ge 0}$. Let $w=w_{-k}+w'$ with $k>0$, $w_{-k}\in D_{-k}\setminus\{0\}$ and $w'\in D_{> -k}$. By $[z,w]=1$ one has $[z_{-s}, w_{-k}]=0$, thus $w_{-k}\in C(z_{-s})=K[z_{-s}]$. Then $s|k$ and $w_{-k}=\al z_{-s}^d$, where $d=k/s$ and $\al\ne 0$. Let $w^{(1)}=w-\al z^d$ and $\mathbf{v}_{-1,1}(w^{(1)})=m$. Then $[z,w^{(1)}]=1$ and $m>-k$. If $m<0$, we continue this procedure, until we get some $w^{(n)}=w-f(z)\in D_{\ge 0}$, where $f\in K[X]$ and $n\ge 1$. It is clear that  $[z,w^{(n)}]=1$.

As $w^{(n)}\in D_{\ge 0}$ and $[z,w^{(n)}]=1$, $z$ and $w^{(n)}$ generate $A_1$ by Theorem \ref{33}. So $z$ and $w$ also generate $A_1$.
\ep
\bco\lb{98}
Assume that $z,w\in A_1$ with $[z,w]=1$.
 If $z\in D_{\ge -1}$, then $z$ and $w$ generate $A_1$.
\eco
\bp
By Theorem \ref{33}, we only need to prove it in the case \[z=z_{-1}+z',\ z_{-1}\in D_{-1}\setminus\{0\},\ z'\in D_{\ge 0}.\] But now $C(z_{-1})=K[z_{-1}]$ by Corollary \ref{51}, so the result follows from the above corollary.
\ep
\begin{prop}\lb{93}
 Assume that $z,w\in A_1$, $[z,w]=1$, and $\{\mathbf{f}_{\rho,\sigma}(z), \mathbf{f}_{\rho,\sigma}(w)\}=1$ for some $(\rho,\sigma)\in \mathbb{R}^2$ with $\rho+\sigma>0$. Then $z$ and $w$ generate $A_1$.
\end{prop}
\bp
Let \[G_1\rt \widetilde{G}_1, \psi\mapsto \tilde{\psi}\] be the isomorphism in \eqref{96}.

Set \[f=\mathbf{f}_{\rho,\sigma}(z), g=\mathbf{f}_{\rho,\sigma}(w).\] Then  $\{f,g\}=1$.

 Assume that $(z^*,w^*)=(\psi(z),\psi(w))$ for some $\psi\in G_1$. Then $[z^*,w^*]=1$.  Let
\[f^*= \tilde{\psi}(f), g^*= \tilde{\psi}(g). \]  Then $\{f^*,g^*\}=1$.

Let \[(\rho^*,\sigma^*)=\left\{
\begin{array}{ll}
(\rho,\sigma) & \text { if } \psi\in G_0; \\
(\sigma,\rho) & \text { if } \psi\in G_0\psi_0.
\end{array}
\right.\]  Then \[f^*=\mathbf{f}_{\rho^*,\sigma^*}(z^*), g^*=\mathbf{f}_{\rho^*,\sigma^*}(w^*).\] One has the following commutative diagram:
\[\begin{CD}
  (z,w) @>\mathbf{f}_{\rho,\sigma} >> (f,g)\\
 @V \psi VV  @V \tilde{\psi}  V  V\\
(z^*,w^*) @>\mathbf{f}_{\rho^*,\sigma^*}>> (f^*,g^*).
\end{CD}\]

Recall that $\eta(z,w)=(w,-z)$. By $\{\mathbf{f}_{\rho,\sigma}(z), \mathbf{f}_{\rho,\sigma}(w)\}=1$, one has $\{\mathbf{f}_{\rho,\sigma}(w), \mathbf{f}_{\rho,\sigma}(-z)\}=1$.

It is clear that $z, w$ generate $A_1$ if and only if $z^*, w^*$ generate $A_1$, if and only if $w,-z$ generate $A_1$. As $U$ is generated by $G_1$ and $\eta$, we only need to prove the result for suitable representative in the $U$-orbit of $(z,w)$, which corresponds to suitable representative in the $\widetilde{U}$-orbit of $(f,g)$.

%We will assume that $z,w$ both have constant term 0.

By Proposition \ref{30}, up to the action of $\widetilde{U}$, $(f,g)$ will be in one of the following 4 cases.

(1) $(f,g)=(X,Y)$.

Assume $\rho\ge 0,\sigma\ge 0, \rho+\sigma>0$. If $\rho\le \sigma$ then $\mathbf{v}_{0,1}(z)\le 1$; If $\sigma\le \rho$ then $\mathbf{v}_{1,0}(w)\le 1$. By Proposition \ref{3}
 $z, w$ generate $A_1$.

%If $\rho=0,\sigma>0$, then $z=p, w=q+h(p)$ for some $h\in K[X]$ and $z, w$ generate $A_1$.

Assume $\rho<0,\rho+\sigma>0$. Then $z\in D_{\le 0}$, and $z, w$ generate $A_1$ by Theorem \ref{33}. The case $\sigma<0,\rho+\sigma>0$ are treated similarly.

(2) $(f,g)=(\alpha X+\beta Y, \gamma X+\delta Y), \alpha,\delta,\beta,\gamma\in K^*, \alpha\delta-\beta\gamma=1$.

Now $(\rho,\sigma)=(1,1)$, $z=\alpha p+\beta q+\lambda, w=\gamma p+\delta q+\mu$ and $z, w$ generate $A_1$.

(3) $(f,g)=(X+\lambda Y^n, Y), \la\ne 0, and\ n\ge 1$.

Now $(\rho,\sigma)=(n,1)$, $z=p+\la q^n+h(q)$ for some $h\in K[X]$ with $deg(h)<n$, and $w=q+\mu$. Then $z, w$ generate $A_1$.

(4) $(f,g)=(X+\al,Y), \al\ne 0$.

Now $(\rho,\sigma)=(0,1)$, $z= p+\al, w=q+h(p)$ for some  $h\in K[X]$, and $z, w$ generate $A_1$.
\ep

\bthm\lb{56}
 Assume that $z\in A_1\setminus\{0\}$ and $f=\mathbf{f}_{\rho,\sigma}(z)$ for some $(\rho,\sigma)\in \mathbb{R}^2$ with $\rho+\sigma>0$. If $C(f)=K[f]$ and $[z,w]=1$ for some $w\in A_1$, then $z$ and $w$ generate $A_1$.
\ethm
%\begin{rem}
%If $\mathbf{v}_{\rho,\sigma}(z)>0$, $C(f)=K[f]$, and $[z,w]=1$ for some $w\in A_1$, then $z$ and $w$ generate $A_1$.
%\end{rem}
\bp
 As $C(f)=K[f]$, by Corollary \ref{90}, $f$ is not a nonzero multiple of some proper power of a nonconstant polynomial.

If $\rho,\sigma$ are linearly dependent over $\mathbb{Q}$, then replace $(\rho,\sigma)$ by some positive multiple  of itself so that  $(\rho,\sigma) \in \mathbb{Z}^2$ with $\rho+\sigma>0$, and $\mathbf{f}_{\rho,\sigma}(z)$ is still $f$.

If $\rho,\sigma$ are linearly independent over $\mathbb{Q}$, then $f$ is a monomial. We can find some $(\rho',\sigma')$ in a sufficiently small neighbourhood of $(\rho,\sigma)$ in $\mathbb{R}^2$ such that $\rho',\sigma'$ are linearly dependent over $\mathbb{Q}$, $\mathbf{f}_{\rho',\sigma'}(z)=f$, and $\rho'+\sigma'>0$. Replace $(\rho,\sigma)$ by a suitable positive multiple of $(\rho',\sigma')$ so that  $(\rho,\sigma) \in \mathbb{Z}^2$ with $\rho+\sigma>0$,  and $\mathbf{f}_{\rho,\sigma}(z)$ is still $f$.

 After the adjustment as above (if necessary), one can always assume that $(\rho,\sigma)\in \mathbb{Z}^2$ with $\rho+\sigma>0$,  $\mathbf{f}_{\rho,\sigma}(z)=f$ is not a nonzero multiple of some proper power of a nonconstant polynomial.

 By assumption, $[z,w]=1$ for some $w\in A_1$. Set $ a=\mathbf{v}_{\rho,\sigma}(z)\in \mathbb{Z}$.

 We now prove that if $a\le 0$ then $z$ and $w$ generate $A_1$. Assume that $a\le 0$. If $\rho,\sigma>0$, then $z$ is a constant, which contradicts to $[z,w]=1$; if $\rho=0,\sigma>0$, then $\mathbf{v}_{0,\sigma}(z)=0$ and $z, w$ generate $A_1$ by Proposition \ref{3}; if $\rho<0,\rho+\sigma>0$, then $a\le 0$ implies that $z\in D_{\le 0}$, thus $z$ and $w$ generate $A_1$ by Theorem \ref{33}. The proof of the remaining cases is similar to the above cases and is omitted. We will assume that $a>0$ from now on.

Set $g=\mathbf{f}_{\rho,\sigma}(w)$ and $b=\mathbf{v}_{\rho,\sigma}(w)$. As $[z,w]=1$, $\{f,g\}=0$ or 1. We start the following process.\\[2mm]

Step 1: Judge whether $b\le 0$.

 If $b\le 0$ then for the same reason as above $z$ and $w$ generate $A_1$, and stop the process;

 If $b>0$, go to Step 2.\\[2mm]

Step 2: Judge whether $\{f,g\}=1$.

If $\{f,g\}=1$, then by Proposition \ref{93}, $z$ and $w$ generate $A_1$. Stop the process.

If $\{f,g\}=0$, go to the next step.\\[2mm]

Step 3: As $\{f,g\}=0$, one has $f^b=\la g^a, \la\ne 0$. Let $d=gcd(a,b), a_0=a/d, b_0=b/d$. By Proposition \ref{27},  there exists some $h\in K[X,Y], f=\gamma h^{a_0}, g=\mu h^{b_0}$, $\gamma\mu\ne 0$. By hypothesis, one must have $a_0=1$, thus $g=\beta f^{b_0}$ for some $\beta\ne 0$. Let $w_1=w-\beta z^{b_0}$. Then $[z,w_1]=1$, and, $z,w$ generate $A_1$ if and only if $z,w_1$ generate $A_1$. It is clear that $\mathbf{v}_{\rho,\sigma}(w_1)<\mathbf{v}_{\rho,\sigma}(w)$.

Then go back to Step 1 and repeat the process for $w_1$. If the process does not stop at Step 1 and Step 2, then we get some $w_2$ at Step 3. Repeat this process and we get $w_1, w_2,\cdots, w_i,\cdots$. As $\mathbf{v}_{\rho,\sigma}(w_{i+1})<\mathbf{v}_{\rho,\sigma}(w_i)$ for all $i$, the process will terminate at Step 1 or Step 2 after finite many steps.

So the proof is concluded.

\ep

\bco\lb{62}  Assume that $z\in A_1\setminus\{0\}$ and $f=\mathbf{f}_{\rho,\sigma}(z)$ for some $(\rho,\sigma)\in \mathbb{R}^2$ with $\rho+\sigma>0$.

 (1) If $|E(f)|=2$, then $z$ and $w$ generate $A_1$.

 (2) Assume that $|E(f)|=1$ and $f=\la X^iY^j, \la\ne 0$. If $i\ge 1, j\ge 1$, and $gcd(i,j)=1$, then $z$ and $w$ generate $A_1$.
\eco
\bp
(1) We show that if $|E(f)|=2$, then $f$ is not a nonzero multiple of some proper power of a polynomial, thus $C(f)=K[f]$ and the result follows from the above theorem. Assume the contrary that $f=\la h^m$, where $\la\ne 0, h\in K[X,Y],$ and $m$ is an integer $\ge 2$. If $h$ is a monomial then $|E(f)|=1$; if $h$ has at least 2 terms, then $h^m$ with $m\ge 2$ will have at least 3 terms, which also contradicts to $|E(f)|=2$.

(2) It is clear that in this case $f$ cannot be written as $\la h^m$, where $\la\ne 0, h\in K[X,Y],$ and $m$ is an integer $>1$. So the result follows from the above theorem.

\ep

The main result of \cite{bl} says that if both $z$ and $w$  are sums of not more than 2 homogeneous elements, then $z$ and $w$ generate $A_1$. One can generalize it as follows.
\bthm\lb{97}
Assume that $z,w\in A_1$ and $[z,w]=1$. Assume that $z$ is a sum of not more than 2 homogeneous elements of $A_1$, then $z$ and $w$ generate $A_1$.
\ethm
\bp
If  $z$ is a homogeneous element, then the result follows from Lemma \ref{2}. Then we only need to consider the case that  $z$ is a sum of 2 homogeneous elements. Write $z=z_1+z_2$ with $z_1,z_2$ homogeneous. It is clear that there exists some unique $(r,s)\in \mathbb{Z}^2$ with $r+s>0$ and $gcd(r,s)=1$, such that $f=\mathbf{f}_{r,s}(z)=\mathbf{f}_{r,s}(z_1)+\mathbf{f}_{r,s}(z_2)$ satisfies $|E(f)|=2$. See Figure 5. So the result follows from Corollary \ref{62} (1).
\ep

\begin{figure} [hb]
  \centering
  \begin{tikzpicture} [scale=0.7]
    \tikzstyle{every node}=[scale=0.7]
    \draw [->] (0,0) -- (6,0) node [below] {$x$};
    \draw [->] (0,0) -- (0,6) node [left] {$y$};

    \coordinate (w1) at (3,0);
    \coordinate (w2) at (4,1);
    \coordinate (z1) at (0,1);
    \coordinate (z2) at (1,2);
    \coordinate (x1) at (4,5);
    \coordinate (x2) at (3,2);

    \draw (w1) -- node [below] {$z_2$} (w2)
    (z1) -- node [ above] {$z_1$} (z2);
    \draw [<-] (x1) -- node [above left] {$ \left(r,s\right)$ } (x2);
    \foreach \j in {1,2}
    {
      \foreach \k in {w,z}
      {
        \fill (\k\j) circle (1.5pt);
      }
    }
    \fill (x2) circle (1.5pt);
    \fill (0,0) circle (1.5pt);
    \draw [dashed] (0,2.33)--(7,0);

    \node [below] at (0,0) {$O$};
  \end{tikzpicture}
  \caption*{Figure 5}
  \end{figure}

\section*{Acknowledgements}

We would like to heartily thank Chengbo Wang for his help during the proof of the main result.

\end{document}